\documentclass[twocolumn,12pt]{article}			 
\newtheorem{theorem}{Theorem}

\newcommand{\ds}{\displaystyle}
\newcommand{\Interi}{\hbox{\bf Z}}

\newcommand{\qed}{\hfill\rule{2mm}{2mm}}  
\newenvironment{proof}{\begin{trivlist}
\item[\hspace{\labelsep}{\bf\noindent Proof. }]
}{\qed\end{trivlist}}
\newenvironment{longabstract}{\begin{trivlist}
\setlength{\linewidth}{5.5in} 
\item[\hspace{\labelsep}{\bf\noindent Abstract: }]
}{\end{trivlist}}
\title{\Large\bf 
		On some transformations of bilateral \\
		birth--and--death processes with applications\\ 
  to first passage time evaluations\thanks{Paper appeared in: \newline
  The 17th Symposium of Information Theory and Its Applications (SITA '94), 
  Hiroshima, Japan, December 6-9, 1994, pp.\ 739--742.}
		}
\author{\normalsize
	{\sc Antonio Di Crescenzo\footnote{current address: \newline 
 Dipartimento di Matematica e Informatica, 
 Universit\`a di Salerno, Via Ponte don Melillo,  I-84084 Fi\-sciano (SA), 
 Italy, Email: adicrescenzo@unisa.it}}
}
\date{}
\begin{document}
\thispagestyle{empty}
\maketitle
\thispagestyle{empty}
\begin{longabstract}
A method yielding simple relationships among bilateral 
birth--and--death processes is outlined. This allows one to relate  
birth and death rates of two processes in such a way 
that their transition probabilities, first--passage--time 
densities and ultimate crossing probabilities are mutually related 
by some product--form expressions. 
\end{longabstract}
\setlength{\baselineskip}{15pt}
%
\section*{\large\bf 1. Introduction}
Birth--and--death processes are a powerful tool to describe stochastic 
models arising in population dynamics (see for instance 
Ricciardi~\cite{Ri}) and also in a large variety of applied fields, 
such as adaptive queueing systems and neurophysiology. 
In particular, their transition probabilities and first--passage--time 
densities play a relevant role in many applied contexts. Hence, 
obtaining closed form expressions for these functions is an important 
task. Unfortunately, apart from a few cases, such expressions are rarely 
encountered in the literature. Furthermore, procedures which have been 
successfully exploited for continuous--state processes are not suitable 
for point processes due to the discreteness of their state--space. 
\par
The aim of this note is to give an extension of a method leading to 
direct transformations between birth--and--death processes 
(cf. Di~Crescenzo~\cite{Di94a} or~\cite{Di94b}) to the 
doubly infinite state--space case. Here we consider a preassigned  
birth--and--death process $X_t$ whose state--space is 
the set of integers. In Section~2 we show how to construct a new 
birth--and--death process $\widetilde X_t$ in such a way that the rates 
of the two processes are mutually related. As a consequence, their 
transition probabilities are related by a simple product--from relation. 
In Section~3 it is shown that also the first--passage--time densities 
of the two processes are similarly connected. 
We point out that the transformation from $X_t$ to $\widetilde X_t$ can 
be viewed as a method to construct new stochastic models. Indeed, making 
use of such method, starting from a birth--and--death process with constant 
rates, in Section~4 a new process is obtained. 
\section*{\large\bf 2. Transition probabilities}
Let $\{X_t;\;t\geq 0\}$ be a birth--and--death process whose state--space 
is the set of integers $\Interi\equiv\{\ldots,-1,0,1,\ldots\}$. 
As usual, for all $n\in \Interi$ we denote by $\lambda_n$ and $\mu_n$ the 
birth and the death rates of $X_t$, i.e.
\begin{eqnarray*}
	&&\lambda_n = \lim_{\delta\downarrow 0}{1\over\delta}\,
	\Pr\{X_{t+\delta}=n+1\,|\,X_t=n\},		\\
	&&\mu_n = \lim_{\delta\downarrow 0}{1\over\delta}\,
	\Pr\{X_{t+\delta}=n-1\,|\,X_t=n\}		
\end{eqnarray*}
for all $n\in\Interi$. We also assume that rates $\lambda_n$ and 
$\mu_n$ are positive, so that the birth--and--death process has no 
absorbing or reflecting states. 
In this case $X_t$ is said to be a {\em bilateral} birth--and--death 
process (cf.~Ismail {\em et~al.}~\cite{IsLeMaVa}). We also assume that $X_t$ 
is {\em simple}; hence, the set of rates $\{\lambda_n,\mu_n\}$ uniquely 
determines the birth--and--death process. According to Callaert and Keilson 
(cf.~\cite{CaKe2}), a bilateral birth--and--death process is 
simple if and only if the two component birth--and--death processes, 
obtained by locating at $n=0$ a boundary reflecting in both directions, 
are simple (necessary and sufficient conditions are also given 
in Pruitt~\cite{Pr}). 
\par
As $X_t$ is simple, the transition probabilities 
$p_{k,n}(t)=\Pr\big\{X_t=n\,|\,X_0=k\big\}$ 
are the unique solution of the following system: 
\begin{eqnarray}
	{d\over dt}p_{k,n}(t)&=&\lambda_{n-1}\,p_{k,n-1}(t)		\label{equazione} \\
	&-& \big(\lambda_n+\mu_n\big)\,p_{k,n}(t)		\nonumber 	\\
	&+& \mu_{n+1}\,p_{k,n+1}(t) \qquad (n\in \Interi). \nonumber
\end{eqnarray}
This has to be solved with initial conditions
\begin{equation}
	p_{k,n}(0)=\cases{
	1 & if $n=k$\cr
	0	& if $n\neq k$,\cr}
	\label{initial}
\end{equation}
where $k\in\Interi$ is the initial state. 
\par
We shall now see that, under suitable assumptions, there exists a bilateral 
birth--and--death process $\widetilde X_t$ having state--space $\Interi$ 
whose rates are obtained from those of $X_t$ and such that the transition 
probabilities of the two processes are mutually related by a product--form 
relation. 
\begin{theorem}\label{tsimilproc}
Let $\{\ldots,\nu_{-1},\nu_0,\nu_1,\ldots\}$ be a strictly monotonic 
sequence of positive numbers satisfying 
\begin{equation}
	\nu_{n+1}\lambda_n-\nu_n(\lambda_n+\mu_n)+\nu_{n-1}\mu_n=0
	\label{condizione}
\end{equation}
for all $n\in \Interi$.	If a bilateral birth--and--death process 
$\widetilde X_t$ has rates  
\begin{equation}
	\widetilde\lambda_n=\lambda_n\,{\nu_{n+1}\over\nu_n}\quad 
	\widetilde\mu_n=\mu_n\,{\nu_{n-1}\over\nu_n}\qquad (n\in \Interi) 
	\label{rates}
\end{equation}
then its transition probabilities are given by:
\begin{equation}
	\widetilde p_{k,n}(t)={\nu_n\over\nu_k}\,p_{k,n}(t)\qquad (k,n\in \Interi). 
	\label{tildeprobtrans}
\end{equation}
\end{theorem}
\begin{proof} 
Let
	$$p^*_{k,n}(t)={\nu_n\over\nu_k}\,p_{k,n}(t)\qquad (k,n\in \Interi).$$ 
Making use of~(\ref{equazione}), for all $n\in\Interi$ we have 
\begin{eqnarray*}
	{d\over dt} p^*_{k,n}(t)&=&{\nu_n\over\nu_k}\,{d\over dt}p_{k,n}(t)	\\
	&=& \lambda_{n-1}\,{\nu_n\over\nu_{n-1}}\,p^*_{k,n-1}(t)	\\
	&-& (\lambda_n+\mu_n)\,p^*_{k,n}(t)	\\
	&+& \mu_{n+1}\,{\nu_n\over\nu_{n+1}}\,p^*_{k,n+1}(t).
\end{eqnarray*}
From relations~(\ref{condizione}) and~(\ref{rates}) it follows  
\begin{eqnarray*}
	\lambda_n+\mu_n &=& \lambda_n\,{\nu_{n+1}\over\nu_n}+\mu_n\,
	{\nu_{n-1}\over\nu_n}	\\
	&=& \widetilde\lambda_n+\widetilde\mu_n\,,	\\
	\lambda_{n-1}\,{\nu_n\over\nu_{n-1}} &=& \widetilde\lambda_{n-1}\,,	\\
	\mu_{n+1}\,{\nu_n\over\nu_{n+1}} &=& \widetilde\mu_{n+1}\,.
\end{eqnarray*}
Hence, the previous system becomes: 
\begin{eqnarray*}
	{d\over dt}p^*_{k,n}(t)&=&\widetilde\lambda_{n-1}\,p^*_{k,n-1}(t)	\\
	&-&(\widetilde\lambda_n+\widetilde\mu_n)\,p^*_{k,n}(t)	\\
	&+&\widetilde\mu_{n+1}\,p^*_{k,n+1}(t)\qquad (n\in\Interi).
\end{eqnarray*}
These are the forward equations of a bilateral birth--and--death process 
$\widetilde X_t$, characterized by rates $\widetilde\lambda_n,\widetilde\mu_n$. 
The corresponding initial conditions are: 
\begin{eqnarray*}
	\lim_{t\downarrow 0}p^*_{k,n}(t)&=&\lim_{t\downarrow 0}{\nu_n\over\nu_k}
	\,p_{k,n}(t)		\\
	&=&{\nu_n\over\nu_k}\,\delta_{k\,n}=\delta_{k\,n}\,.
\end{eqnarray*}
Due to the uniqueness of the solution of the forward equations,  
probabilities $p^*_{k,n}(t)$ and $\widetilde p_{k,n}(t)$ 
coincide, so that~(\ref{tildeprobtrans}) holds. 
\end{proof}
\par
It should be noted that if the ratio $\lambda_n/\mu_n$ is a constant, 
a solution of system~(\ref{condizione}) can be easily obtained. Indeed, 
if there exists a positive constant $c$ $(c\neq 1)$ such that 
$\mu_n=c\,\lambda_n$ for all $n\in\Interi$, then one has: 
	$$\nu_n=1+\beta\,c^n\qquad (n\in\Interi),$$
with $\beta>0$. 
%
\section*{\large\bf 3. First--passage--time densities}
Let 
\begin{eqnarray*}
	&& g_{k,s}(t)=		\\
	&& ={d\over dt}\Pr\big(\inf\{t>0: X_t=s\}<t\,|\,X_0=k\big).
\end{eqnarray*}
be the first--passage--time probability density function 
of $X_t$ through the state $s$, conditioned upon the initial condition 
$\Pr\{X_0=k\}=1$, with $k\neq s$. 
In the following Theorem we determine the relations existing between 
the first--passage--time densities of the birth--and--death processes 
$X_t$ and $\widetilde X_t$. 
\begin{theorem}\label{teoremag1}
Under the assumptions of Theorem~$\ref{tsimilproc}$, the first--passage--time 
density $\widetilde g_{k,s}(t)$ of process $\widetilde X_t$ with 
rates~$(\ref{rates})$ is given by
\begin{equation}
	\widetilde g_{k,s}(t)={\nu_s\over\nu_k}\,g_{k,s}(t)\qquad (k\neq s).
	\label{relazioneg_1}
\end{equation}
\end{theorem}
\begin{proof}
The transition probabilities and the first--passage--time density 
of $\widetilde X_t$ are related by the integral equation 
	$$\widetilde p_{k,n}(t)=\int_0^t \widetilde g_{k,s}(\vartheta)\,
	\widetilde p_{s,n}(t-\vartheta)\,d\vartheta,$$
which holds as $k<s\leq n$ or $n\leq s<k$. 
Making use of relation~(\ref{tildeprobtrans}) one obtains: 
	$$p_{k,n}(t)={\nu_k\over\nu_s}\int_0^t \widetilde g_{k,s}(\vartheta)\,
	p_{s,n}(t-\vartheta)\,d\vartheta.$$
However, one also has: 
	$$p_{k,n}(t)=\int_0^t g_{k,s}(\vartheta)\,
	p_{s,n}(t-\vartheta)\,d\vartheta.$$
Comparing the last two integral equations, due to the uniqueness of their 
solution, identity~(\ref{relazioneg_1}) follows. 
\end{proof}
Bearing in mind that a first--passage--time density could be ``defective'', 
let us now pinpoint some further relations. If sequence $\{\nu_n\}_n$ is 
increasing, from~(\ref{relazioneg_1}) one has: 
\begin{itemize}
 	\item $\widetilde g_{k,s}(t)>g_{k,s}(t)\;\;$ if $\;s>k$; 
 	\item $\widetilde g_{k,s}(t)<g_{k,s}(t)\;\;$  if $\;s<k$. 
\end{itemize}
Instead, when $\{\nu_n\}_n$ is decreasing the inequalities between the 
first--passage--time densities are inverted. 
\par
For all $k\neq s$ let us now introduce the ultimate crossing probabilities
	$${\cal P}_{k,s}=\int_0^{\infty} g_{k,s}(t)\,dt,\quad
	\widetilde {\cal P}_{k,s}=\int_0^{\infty} \widetilde g_{k,s}(t)\,dt.$$
Integration of both sides in~(\ref{relazioneg_1}) over $(0,\infty)$ 
immediately yields
	$$\widetilde {\cal P}_{k,s}={\nu_s\over\nu_k}\,{\cal P}_{k,s}
	\qquad (k\neq s).$$
As $\{\nu_n\}_n$ is a strictly monotonic sequence, the ultimate crossing 
probabilities ${\cal P}_{k,s}$ and $\widetilde {\cal P}_{k,s}$ are 
never equal; hence, in particular, they cannot be simultaneously unity. 
%
\section*{\large\bf 4. An example}
In this Section we give an example of application of previous results. 
Let $X_t$ be the bilateral birth--and--death process with constant rates  
$\lambda_n=\lambda>0$, $\mu_n=\mu>0$ $(n\in\Interi)$. 
The transition probabilities and the first--passage--time density of $X_t$ 
are given by (cf.~Abate {\em et al.}~\cite{AbKiWh}): 
\begin{equation}
	p_{k,n}(t)=e^{-(\lambda+\mu)t}\Bigl({\lambda\over\mu}\Bigr)^{(n-k)/2}\,
	I_{n-k}\bigl(2t\sqrt{\lambda\mu}\bigr)
	\label{probAesempio}
\end{equation}
and
\begin{eqnarray}
	&& g_{k,s}(t)={|s-k|\over t}\,p_{k,s}(t)	\label{probBesempio}		\\
	&& ={|s-k|\over t}\,e^{-(\lambda+\mu)t}\Bigl({\lambda\over\mu}
	\Bigr)^{(s-k)/2}\,I_{s-k}\big(2t\sqrt{\lambda\mu}\big),	\nonumber	
\end{eqnarray}
respectively, where $I_k$ is the modified Bessel function of first kind: 
	$$I_k(2z)=\sum_{j\,=\,0}^{\infty}\,{z^{\,k+2j}\over j!\,(k+j)!}\,.$$
Note that when $\lambda=\mu$ equation~(\ref{condizione}) does not admit 
of a positive non--constant solution. Instead, if $\lambda\neq\mu$ 
Theorem~\ref{tsimilproc} holds. Furthermore, from equality~(\ref{condizione}) 
for all $n\in \Interi$ one has:  
	$$\nu_n=1+\beta\,c^n\qquad(\beta>0),$$
where $c=\mu/\lambda$. Recalling expressions~\ref{rates}, for all 
$n\in\Interi$ the rates of the birth--and--death process $\widetilde X_t$ 
are given by: 
	$$\widetilde\lambda_n=\lambda\,\ds{1+\beta\,c^{n+1}
	\over 1+\beta\,c^{n}}\,,
	\qquad\widetilde\mu_n=\mu\,
	\ds{1+\beta\,c^{n-1}\over 1+\beta\,c^{n}}.$$
Furthermore, making use of equations~(\ref{tildeprobtrans}) 
and~(\ref{probAesempio}), the transition probabilities of 
$\widetilde X_t$ follow: 
\begin{eqnarray*}
	&&\widetilde p_{k,n}(t)=\ds{1+\beta\,c^n\over 1+\beta\,c^k}\,	\\
	&&\qquad\qquad \times 
	e^{-(\lambda+\mu)t}c^{(k-n)/2}\,I_{n-k}\big(2t\sqrt{\lambda\mu}\big).
\end{eqnarray*}
Finally, making use of Theorem~\ref{teoremag1}, 
from~(\ref{probBesempio}) we obtain the first--passage--time 
density of $\widetilde X_t$ from state $k$ to state $s$: 
\begin{eqnarray*}
	&&\widetilde g_{k,s}(t)
	=\ds{1+\beta\,c^s\over 1+\beta\,c^k}\,{|s-k|\over t}	\\
	&&\qquad\qquad \times e^{-(\lambda+\mu)t}c^{(k-s)/2}\,
	I_{s-k}\big(2t\sqrt{\lambda\mu}\big).
\end{eqnarray*}
%

\vfill\eject
\end{document}